\documentclass[12pt]{article}
\usepackage{amssymb}
\usepackage{amsmath}
\usepackage{amsthm}
\usepackage{authblk}
\usepackage{color}
\usepackage{comment}
\usepackage{geometry}		
\usepackage{graphicx}
\usepackage[sort&compress]{natbib}     

\makeatletter
	\@addtoreset{equation}{section}
\makeatother

\let\originalleft\left
\let\originalright\right
\renewcommand{\left}{\mathopen{}\mathclose\bgroup\originalleft}
\renewcommand{\right}{\aftergroup\egroup\originalright}

\title{Grazing bifurcations of linear impact oscillators in the zero damping limit.}
\author{O.J.~Goodman}
\author{D.J.W.~Simpson}
\affil{School of Mathematical and Computational Sciences, Massey University, Palmerston North, New Zealand}

\begin{document}

\newcommand{\rD}{{\rm D}}
\newcommand{\re}{{\rm e}}

\newtheorem{theorem}{Theorem}[section]
\newtheorem{corollary}[theorem]{Corollary}
\newtheorem{lemma}[theorem]{Lemma}
\newtheorem{proposition}[theorem]{Proposition}

\theoremstyle{definition}
\newtheorem{definition}{Definition}[section]
\newtheorem{assumption}[definition]{Assumption}

\theoremstyle{remark}
\newtheorem{remark}{Remark}[section]

\maketitle





\begin{abstract}

We consider a harmonically forced linear impact oscillator, where impact events are instantaneous with energy loss. We study the dynamics at the grazing bifurcation of the non-impacting periodic solution in the limit that the damping coefficient of the oscillator is zero. Through numerical computations we show that a recurring sequence of bifurcations exists between points of resonance. Specifically, resonance creates a stable periodic solution that subsequently loses stability in a secondary grazing bifurcation, then regains stability in a saddle-node bifurcation, then transitions to a chaotic attractor through a period-doubling cascade. The dynamics persist under mild parameter variation, so apply to weakly-damped impact oscillators near grazing.

\end{abstract}

\section{Introduction}
\label{sec:intro}

Many mechanical systems have rigid vibrating components that occasionally impact each other.
For such systems, recurring impacts arise through {\em grazing bifurcations} \cite{Br99,BlCz99,Ib09}.
These are the critical parameter values at which
a stable periodic motion involves components
that come together for an instant, touching with zero relative velocity before moving apart.
That is, grazing bifurcations delineate a near miss from an impact.
Since impacts are highly nonlinear,
grazing bifucations often vastly increase the complexity
of the motion, such as by generating chaos \cite{DaZh07,MoCh20,HaWi07,DaZh05}.

\begin{figure}[t!]
\begin{center}
\includegraphics[width=10.5cm]{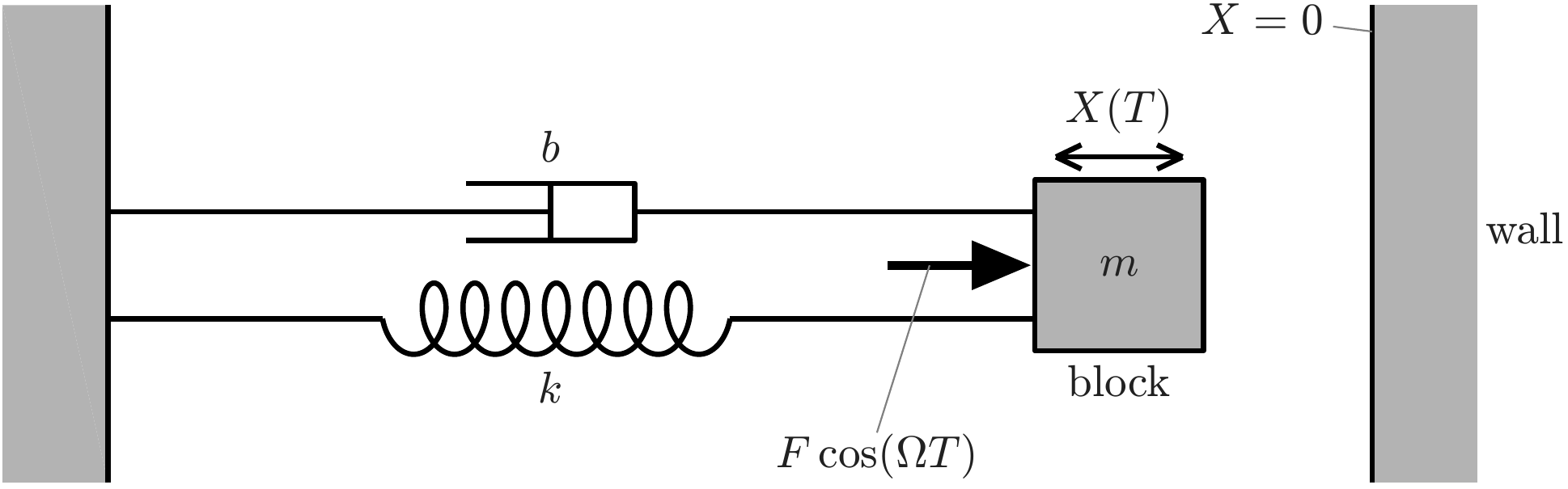}
\caption{
The impact oscillator system modelled by \eqref{eq:modelOriginal}.
The block has position $X(T) \le 0$, mass $m$, damping coefficient $b$, spring constant $k$,
and equilibrium position $X_{\rm eq} < 0$.
Whenever the block hits the wall at $X=0$,
the velocity of the block is reversed and scaled by the coefficient of restitution $\epsilon$.
\label{fig:oscSchem}
} 
\end{center}
\end{figure}

To understand grazing bifurcations in diverse applications,
it is helpful to study minimal systems
that contain only the features necessary
for grazing bifurcations to occur in generic fashion.
A well-studied example of such a system
is the harmonically forced, one-degree-of-freedom,
linear oscillator depicted in Fig.~\ref{fig:oscSchem}.
The equations of motion are
\begin{equation}
\begin{split}
m \frac{d^2 X}{d T^2} + b \frac{d X}{d T} + k (X - X_{\rm eq})
&= F \cos(\Omega T), \qquad \text{while $X < 0$}, \\
\frac{d X}{d T} &\mapsto -\epsilon \frac{d X}{d T}, \qquad \text{whenever $X = 0$},
\end{split}
\label{eq:modelOriginal}
\end{equation}
where impact events are treated as instantaneous
with coefficient of restitution $\epsilon \in [0,1]$.

For sufficiently low values of the forcing amplitude $F$,
the system \eqref{eq:modelOriginal} has a
stable period-$\frac{2 \pi}{\Omega}$ solution with no impacts.
At a critical value of $F$, this solution undergoes a grazing bifurcation by reaching the wall.
The dynamics created in this bifurcation were studied by
several groups in the 1990's \cite{PeVa92,Iv93,Fo94,Pe96}.
They observed that typical grazing bifurcations create chaotic attractors,
and that a stable period-$\frac{2 \pi}{\Omega}$ solution with one impact per period
could only arise when the forcing frequency and damped natural frequency obey a certain resonance.
Later, Nordmark \cite{No01} and others \cite{ThZh06,SiGh26}
found there exist other resonances where 
a stable period-$\frac{2 \pi p}{\Omega}$ solution
is created with some $p \ge 2$, again with one impact per period.
Although these periodic solutions are only created at resonant grazing bifurcations,
they are often stable near typical grazing bifurcations.
Despite the simplicity of \eqref{eq:modelOriginal},
a complete understanding of its dynamics is at present out of reach
as it contains many independent parameters and has a remarkably complex bifurcation structure. 

However, in many physical systems
for which grazing bifurcations have been studied,
the oscillator has comparatively small damping.
Pavlovskaia {\em et al.}~\cite{PaIn10} estimate a dimensionless {\em damping ratio} $\zeta = \frac{b}{2 \sqrt{k m}}$
of $\zeta \approx 0.01$ for their physical experiments.
de Weger {\em et al.}~\cite{DeVa00} also measured $\zeta \approx 0.01$, while
Witelski {\em et al.}~\cite{WiVi14} give $\zeta = 0.0041$ and $\zeta = 0.0043$ for pendula of different sizes.

\begin{figure}[t!]
\begin{center}
\includegraphics[width=\textwidth]{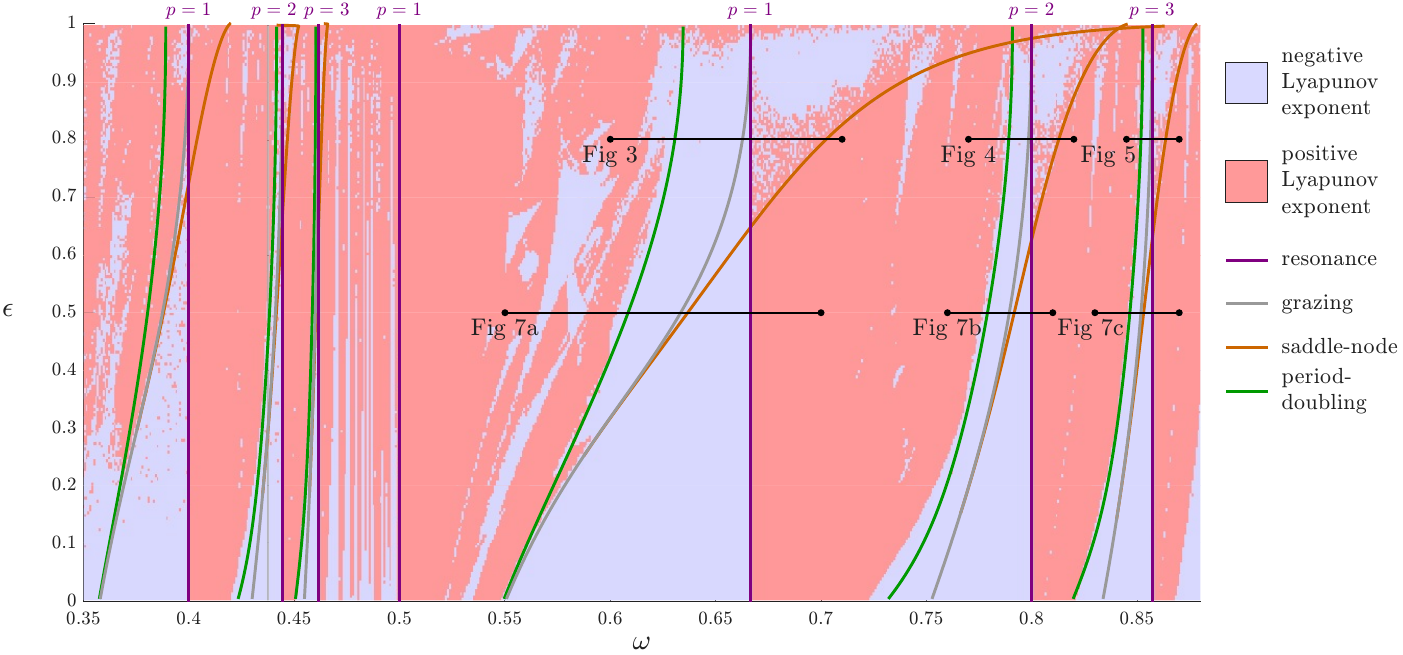}
\caption{
A two-parameter bifurcation diagram of the
impact oscillator system \eqref{eq:modelOriginal} at grazing and in the zero damping limit.
In this case, 
the system has the nondimensionalised form \eqref{eq:modelReduced},
and there are two parameters: the nondimensionalised forcing frequency $\omega$,
and the coefficient of restitution $\epsilon$.
Each point in a $600 \times 200$ grid is coloured
according to the sign of a maximal Lyapunov exponent
obtained numerically from the forward orbit of a random initial condition
(the run-time of this simulation was about $100$ hours).
Numerically computed bifurcation curves of some stable period-$\frac{2 \pi p}{\omega}$ solutions
are overlaid.
One-parameter bifurcation diagrams corresponding to the six horizontal line segments
are discussed in \S\ref{sec:bif}.
\label{fig:bifSet}
} 
\end{center}
\end{figure}

In this paper we study the dynamics of the linear impact oscillator model \eqref{eq:modelOriginal}
at the grazing bifurcation when the damping coefficient of the oscillator is zero.
In this scenario, only two dimensionless parameters remain:
the nondimensionalised forcing frequency $\omega = \frac{\sqrt{m} \Omega}{\sqrt{k}}$,
and the coefficient of restitution $\epsilon$.
Fig.~\ref{fig:bifSet} summarises the stable dynamics
over a range of values of $\omega$ and $\epsilon$.
Most significantly, we find that stable period-$\frac{2 \pi p}{\Omega}$ solutions
are created at points of resonance, then transition to chaos via
a consistent sequence of bifurcations.

Importantly, the dynamics broadly persists
for sufficiently small changes to the parameter values.
This is because the periodic solutions are generically hyperbolic,
and the chaotic attractors may be robust due to the square-root singularity associated
with grazing events \cite{LiCh17,MiLi19}.
Therefore, the dynamics indicated in Fig.~\ref{fig:bifSet}
occurs also for the weakly-damped oscillators near grazing.
In this way, Fig.~\ref{fig:bifSet} explains
how typical grazing bifurcations modify the motion of the linear impact oscillator
when the damping coefficient is small.

The remainder of this paper is organised as follows.
In \S\ref{sec:prelim} we nondimensionalise the equations of motion
and compute the grazing bifurcation of the non-impacting periodic solution.
Here we also describe the resonances that lead to stable periodic impacting motion.

In \S\ref{sec:bif} we study the dynamics
of the system at the grazing bifurcation when the damping is zero.
We employ Monte-Carlo simulations and numerical continuation
to explain the existence and coexistence of stable periodic solutions and chaotic motions,
as summarised by Fig.~\ref{fig:bifSet}.
In \S\ref{sec:lowDamping} we relate the dynamics of Fig.~\ref{fig:bifSet}
with those of the system with low damping near grazing.
Conclusions are presented in \S\ref{sec:conc}.

\section{Preliminaries}
\label{sec:prelim}

\subsection{Nondimensionalisation}

Under the linear scaling of space and time
\begin{align}
x &= \frac{X}{|X_{\rm eq}|}, &
t &= \frac{\sqrt{k} \, T}{\sqrt{m}},
\end{align}
the system \eqref{eq:modelOriginal} becomes
\begin{equation}
\begin{split}
\ddot{x} + 2 \zeta \dot{x} + x + 1 &= A \cos(\omega t), \qquad \text{while $x < 0$}, \\
\dot{x} &\mapsto -\epsilon \dot{x}, \qquad \text{whenever $x = 0$},
\end{split}
\label{eq:modelNondim}
\end{equation}
where dots denote differentiation with respect to $t$, and
\begin{align}
\zeta &= \frac{b}{2 \sqrt{k m}}, &
A &= \frac{F}{k |X_{\rm eq}|}, &
\omega &= \frac{\sqrt{m} \, \Omega}{\sqrt{k}}.
\label{eq:zetaAomega}
\end{align}
The constants \eqref{eq:zetaAomega} are dimensionless;
$\zeta$ is the damping ratio,
while $A$ and $\omega$ are the nondimensionalised forcing amplitude and frequency.
The equilibrium position of the nondimensionalised oscillator is $x = -1$.

\subsection{Grazing}

If impact events are ignored,
then the system with $\zeta > 0$
has a globally attracting period-$\frac{2 \pi}{\omega}$ solution.
The maximum $x$-value attained by this solution is
\begin{equation}
x_{\rm max} = -1 + \frac{A}{\sqrt{\left( 1-\omega^2 \right)^2 + 4 \zeta^2 \omega^2}}.
\nonumber
\end{equation}
Thus if $x_{\rm max} < 0$, then this solution is a non-impacting
solution of the full system \eqref{eq:modelNondim}.
Grazing occurs when $x_{\rm max} = 0$.
Thus
\begin{equation}
A_{\rm graz} = \sqrt{\left( 1-\omega^2 \right)^2 + 4 \zeta^2 \omega^2},
\label{eq:Agraz}
\end{equation}
is the value of $A$ at which
the non-impacting periodic solution undergoes a grazing bifurcation.

\subsection{Return maps and resonance}

To analyse periodic solutions created by the grazing bifurcation $A = A_{\rm graz}$,
one can study a Poincar\'e map or stroboscopic map that captures the local dynamics.
These maps are piecewise-smooth,
where one piece of the map corresponds to an oscillation without an impact,
and the other piece of the map corresponds to an oscillation with an impact.

The piece of the map associated with impacts contains a square-root singularity \cite{No91}.
Thus the derivative of this piece of the map
tends to infinity as the impact velocity tends to zero.
Any local periodic solution created in the grazing bifurcation
involves an impact with vanishingly small impact velocity,
thus is unstable, except in special cases.
These special cases occur when the forcing frequency of \eqref{eq:modelNondim}
is a certain rational multiple of the {\em damped natural frequency}
\begin{equation}
\omega_1 = \sqrt{1 - \zeta^2}.
\label{eq:dnf}
\end{equation}
Ivanov \cite{Iv93} showed that an impacting period-$\frac{2 \pi}{\omega}$ solution
can only be stable as it issues from the grazing bifurcation if
\begin{equation}
\frac{\omega_1}{\omega} = \frac{n}{2},
\label{eq:res1}
\end{equation}
for some $n \in \mathbb{Z}$.
Later Nordmark \cite{No01} studied period-$\frac{2 \pi p}{\omega}$ solutions
with one impact per period, where $p \ge 2$ is an integer.
Such a solution can only be stable as it issues from the grazing bifurcation if
\begin{equation}
\frac{\omega_1}{\omega} = n + \frac{1}{2 p},
\label{eq:res2}
\end{equation}
for some $n \in \mathbb{Z}$, see \cite{SiGh26}.

\section{Dynamics at grazing with zero damping}
\label{sec:bif}

In this section we study the nondimensionalised system \eqref{eq:modelNondim}
with $A = A_{\rm graz}$ and $\zeta = 0$.
In this case, \eqref{eq:modelNondim} reduces to
\begin{equation}
\begin{split}
\ddot{x} + x + 1 &= \left| 1 - \omega^2 \right| \cos(\omega t), \qquad \text{while $x < 0$}, \\
\dot{x} &\mapsto -\epsilon \dot{x}, \qquad \text{whenever $x = 0$}.
\end{split}
\label{eq:modelReduced}
\end{equation}

\subsection{Resonant frequencies}
\label{sub:pointsOfResonance}

With $\zeta = 0$, the natural frequency is $\omega_1 = 1$.
Thus, by \eqref{eq:res1}, $p=1$ resonance occurs when the frequency $\omega$ equals
\begin{equation}
\omega_{1,n} = \frac{2}{n},
\label{eq:omega_1n}
\end{equation}
for some integer $n \ge 1$.
Over the range of $\omega$-values plotted in Fig.~\ref{fig:bifSet},
there are three such frequencies.
These are the values $\omega_{1,3} = \frac{2}{3}$,
$\omega_{1,4} = \frac{1}{2}$, and
$\omega_{1,5} = \frac{2}{5}$,
indicated with purple lines and labelled `$p=1$' in Fig.~\ref{fig:bifSet}.

By \eqref{eq:res2}, $p \ge 2$ resonance occurs at frequencies
\begin{equation}
\omega_{p,n} = \frac{2 p}{1 + 2 n p},
\label{eq:omega_pn}
\end{equation}
for some integer $n \ge 1$.
Over the range of $\omega$-values plotted in Fig.~\ref{fig:bifSet},
we indicate all such frequencies having $p = 2$ or $p = 3$; specifically
$\omega_{2,1} = \frac{4}{5}$,
$\omega_{2,2} = \frac{4}{9}$,
$\omega_{3,1} = \frac{6}{7}$, and
$\omega_{3,2} = \frac{6}{13}$.
The other resonant frequencies present over the plotted range
are $\omega_{p,2}$ with $p \ge 4$, which converge to $\omega = \frac{1}{2}$ as $p \to \infty$.

\subsection{Bifurcation diagrams}
\label{sub:bifDiags}

\begin{figure}[b!]
\begin{center}
\includegraphics[width=\textwidth]{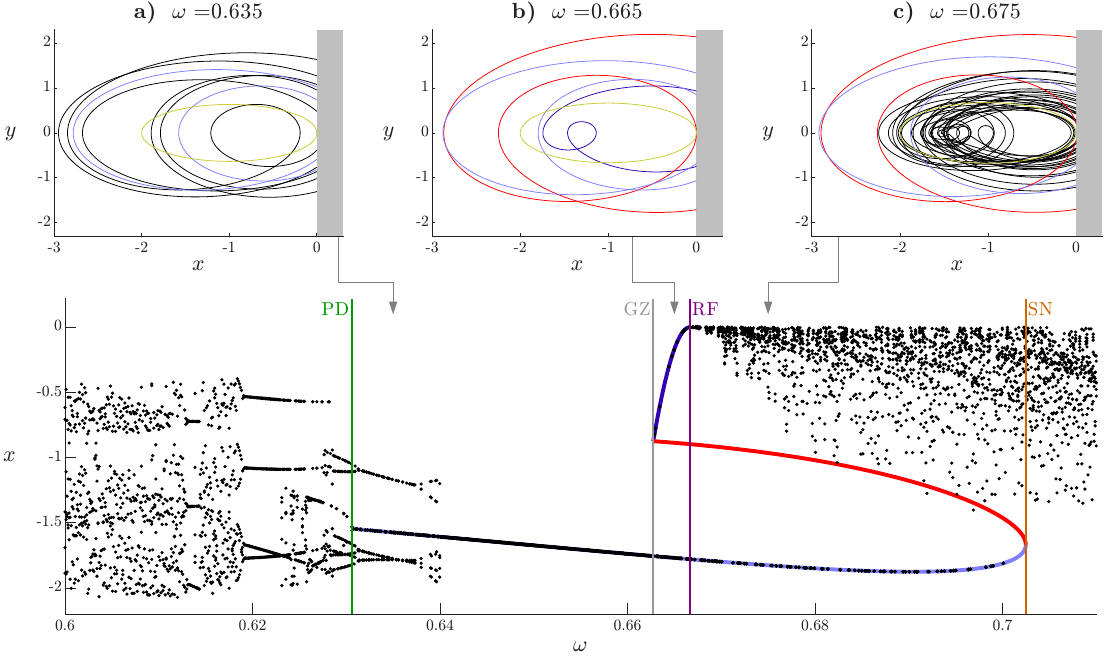}
\caption{
The lower plot is a stroboscopic
bifurcation diagram of \eqref{eq:modelReduced} with $\epsilon = 0.8$
(PD: period-doubling bifurcation;
GZ: grazing bifurcation;
RF: resonant frequency; 
SN: saddle-node bifurcation).
On the vertical axis we have plotted the $x$-values
of orbits at integer multiple times of the period $\frac{2 \pi}{\omega}$ of the forcing.
We have numerically continued the periodic orbit from resonance to period-doubling
(dark blue: stable single-impact;
red: unstable two-impact;
light blue: stable two-impact).
We have also overlaid the result of a Monte-Carlo simulation,
showing stroboscopic points of
orbits of random initial points after transient dynamics has decayed.
The phase portraits show the computed solutions at three different values of $\omega$,
and the grazing non-impacting solution in yellow.
\label{fig:bifDiagPPd}
} 
\end{center}
\end{figure}

\begin{figure}[b!]
\begin{center}
\includegraphics[width=\textwidth]{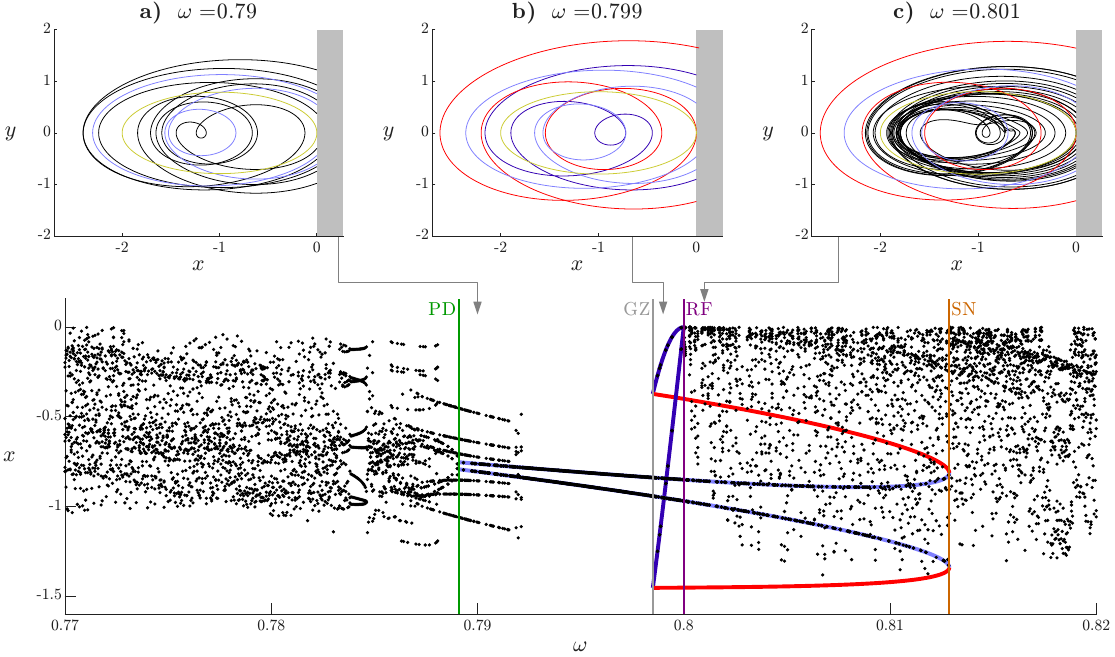}
\caption{
A bifurcation diagram of \eqref{eq:modelReduced} with $\epsilon = 0.8$
using the same conventions as Fig.~\ref{fig:bifDiagPPd}.
The $p=2$ resonant frequency (RF) creates a stable, single-impact,
period-$\frac{4 \pi}{\omega}$ solution (dark blue).
This changes to an unstable two-impact period-$\frac{4 \pi}{\omega}$ solution (red),
then to a stable two-impact period-$\frac{4 \pi}{\omega}$ solution (light blue).
\label{fig:bifDiagPPe}
} 
\end{center}
\end{figure}

\begin{figure}[b!]
\begin{center}
\includegraphics[width=\textwidth]{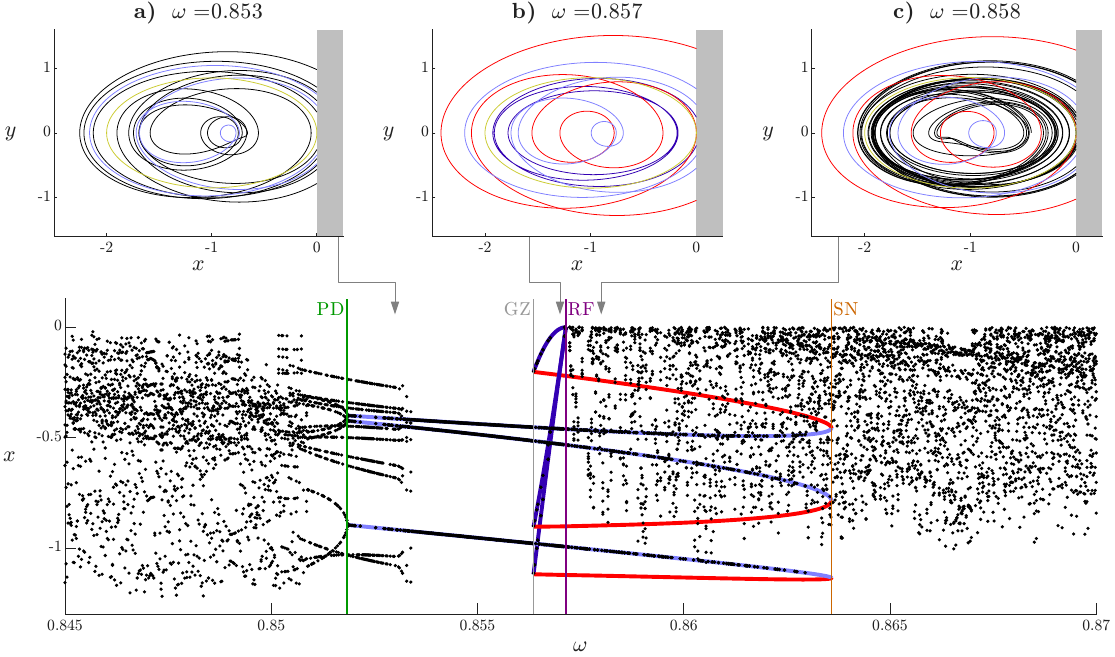}
\caption{
A bifurcation diagram of \eqref{eq:modelReduced} with $\epsilon = 0.8$
using the same conventions as Fig.~\ref{fig:bifDiagPPd}.
The $p=3$ resonant frequency (RF) creates a stable, single-impact,
period-$\frac{6 \pi}{\omega}$ solution (dark blue).
This changes to an unstable two-impact period-$\frac{6 \pi}{\omega}$ solution (red),
then to a stable two-impact period-$\frac{6 \pi}{\omega}$ solution (light blue).
\label{fig:bifDiagPPf}
} 
\end{center}
\end{figure}

The significance of the resonant frequencies
is that a local, stable, single-impact, periodic solution emerges
as the value of $A$ is varied from $A_{\rm graz}$.
In Fig.~\ref{fig:bifSet}, $A = A_{\rm graz}$ is fixed,
yet the resonant frequencies are still important.
For all seven resonant frequencies marked with a purple line,
except $\omega_{1,4} = \frac{1}{2}$
which behaves differently as explained in \cite{SiGh26},
we find that a local, stable, single-impact, period-$\frac{2 \pi p}{\omega}$ solution
emerges as the value of $\omega$ is decreased from the resonant frequency.
This is shown in Figs.~\ref{fig:bifDiagPPd}, \ref{fig:bifDiagPPe}, and \ref{fig:bifDiagPPf}
for the resonant frequencies $\omega_{1,3} = \frac{2}{3}$, $\omega_{2,1} = \frac{4}{5}$, and
$\omega_{3,1} = \frac{6}{7}$ when $\epsilon = 0.8$.
These figures contain bifurcation diagrams
produced by evaluating orbits
{\em stroboscopically}: at integer multiple times of the forcing period $\frac{2 \pi}{\omega}$.
Each figure contains three phase portraits,
and in each figure the bifurcating stable period-$\frac{2 \pi p}{\omega}$ solution
is coloured dark blue in middle phase portrait. 

If the value of $\omega$ is instead increased from the resonant frequency,
then a local chaotic attractor appears to be created.
In Figs.~\ref{fig:bifDiagPPd}--\ref{fig:bifDiagPPf},
this chaotic attractor is coloured black in the right-most phase portrait. 
The phase portraits also show the non-impacting period-$\frac{2 \pi}{\omega}$ solution in yellow.
This solution grazes $x=0$ because $A = A_{\rm graz}$.

These observations explain why in Fig.~\ref{fig:bifSet} the Lyapunov exponent
of numerically computed orbits changes sign along much of the purple lines.
To the left of resonance, there is a stable periodic solution with a negative Lyapunov exponent,
while to the right of resonance, there is a chaotic attractor with a positive Lyapunov exponent.

\subsection{Bifurcations leading to fragile chaos}
\label{sub:bifs}

In Figs.~\ref{fig:bifDiagPPd}--\ref{fig:bifDiagPPf}
we have numerically continued the period-$\frac{2 \pi p}{\omega}$ solution
that issues from resonance through a sequence of bifurcations
leading ultimately to fragile chaos.
Interestingly, the same sequence of bifurcations seems to occur
for every point of resonance (except $\omega = \frac{1}{2}$),
and for all values of $\epsilon$.

When the stable period-$\frac{2 \pi p}{\omega}$ solution 
is created at resonance, it has one impact per period.
As the value of $\omega$ is decreased, this solution,
as a set in $(x,y)$-phase space,
varies position and shape in a continuous fashion.
In doing so, it quickly develops an additional loop
near the equilibrium value $x=-1$
(this loop is seen in Fig.~\ref{fig:bifDiagPPd}b and Fig.~\ref{fig:bifDiagPPe}b),
As $\omega$ is decreased further, the loop expands and soon grazes $x=0$.
This is a grazing bifurcation of nonsmooth-fold type \cite{DiBu08}.
It generates an unstable period-$\frac{2 \pi p}{\omega}$ solution
that has two impacts per period and coexists with the single-impact solution.

As the value of $\omega$ is increased from this grazing bifurcation,
the unstable two-impact solution persists
and is plotted red in Figs.~\ref{fig:bifDiagPPd}--\ref{fig:bifDiagPPf}.
This solution continues until regaining stability in a saddle-node bifurcation.
It then extends left from the saddle-node bifurcation
until losing stability in a period-doubling bifurcation.
While stable, the two-impact period-$\frac{2 \pi p}{\omega}$ solution
is coloured light blue and shown in all three phase portraits in each of
Figs.~\ref{fig:bifDiagPPd}--\ref{fig:bifDiagPPf}.

In each figure, the period-doubling bifurcation
appears to initiate a period-doubling cascade to chaos.
Since the chaos is produced by a period-doubling cascade,
we expect that the chaos is fragile in the sense that chaotic attractors
only exist at isolated values of $\omega$
because periodic windows are dense in parameter space.

\subsection{Robust chaos}
\label{sub:robustChaos}

\begin{figure}[b!]
\begin{center}
\includegraphics[width=\textwidth]{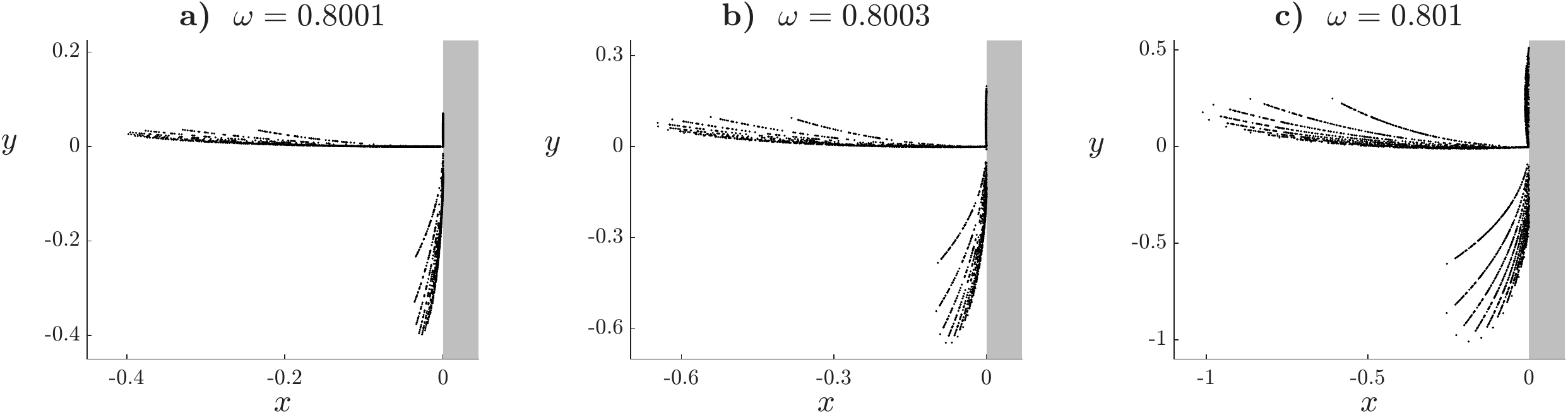}
\caption{
Stroboscopic maps showing the apparently chaotic attractor
of \eqref{eq:modelReduced} with $\epsilon = 0.8$ and three different values of $\omega$. 
\label{fig:robustChaos}
} 
\end{center}
\end{figure}

As mentioned above, as the value of $\omega$ is increased from resonance,
a local chaotic attractor appears to be created.
This attractor is shown in Fig.~\ref{fig:robustChaos}
at three values of $\omega$ just beyond $\omega_{2,1} = \frac{4}{5}$.
Each plot shows the $(x,y)$ values of a
numerically computed orbit evaluated stroboscopically.

In the $(x,y)$-plane, the size of the attractor increases from zero at resonance.
As evidenced in Fig.~\ref{fig:bifDiagPPe}, the size of the attractor
appears to be asymptotically proportional to the square-root
of the change in the value of $\omega$.
This is consistent with the behaviour
of families of one-dimensional maps with a square-root singularity \cite{ChOt94,No97}.
As $\omega$ increases, the shape of the attractor remains broadly the same,
consisting of three limbs that get fatter with $\omega$,
revealing a fractal substructure within each limb.
Similar behaviour is observed immediately beyond other resonant frequencies.

Based on these observations, we expect that the attractor
has no periodic windows to the right of $\omega_{p,n}$.
That is, the chaotic attractor is robust over an interval of $\omega$-values
with left endpoint $\omega_{p,n}$.
Such robustness has been proved
for some families of two-dimensional square-root maps \cite{MiLi19}.
It remains for future work to formally verify the robustness observed here.

\subsection{Less coexistence of attractors at lower values of $\epsilon$}
\label{sub:more}

The numerical results discussed in this section have thus far considered only $\epsilon = 0.8$.
As $\epsilon$ is varied from $0.8$,
the $\omega$-values of the grazing, saddle-node,
and period-doubling bifurcations follow the curves in Fig.~\ref{fig:bifSet}
which were computed by numerical continuation.
Notice that as the value of $\epsilon$ is decreased,
all of the bifurcations shift to the left.
In particular, each saddle-node bifurcation
passes through the associated resonant frequency,
and converges to the corresponding grazing bifurcation as $\epsilon \to 0$.

\begin{figure}[b!]
\begin{center}
\includegraphics[width=\textwidth]{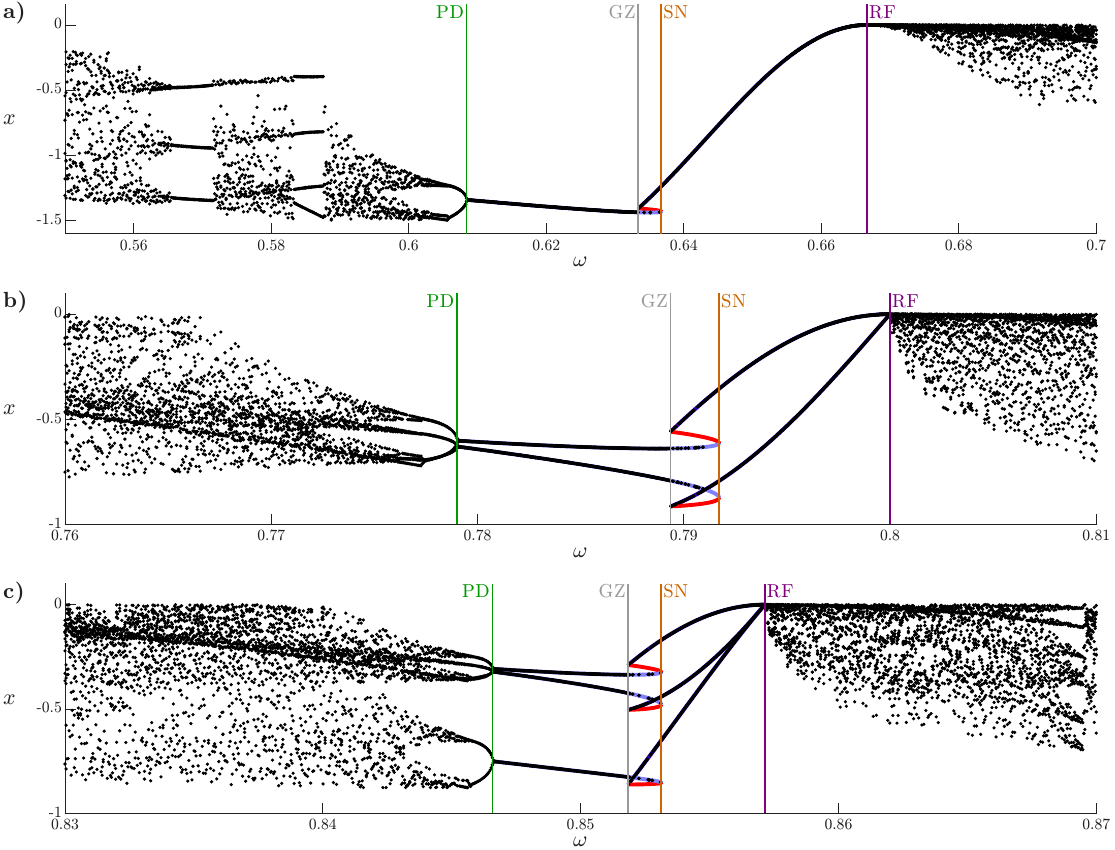}
\caption{
Bifurcation diagrams of \eqref{eq:modelReduced} with $\epsilon = 0.5$.
\label{fig:bifDiagabc}
} 
\end{center}
\end{figure}

Fig.~\ref{fig:bifDiagabc} shows bifurcation diagrams at $\epsilon = 0.5$,
for which the saddle-node bifurcations are already
substantially closer to the grazing bifurcations.
This is significant because the grazing and saddle-node bifurcations
bound an interval of coexistence where the system
has two stable periodic solutions.
There is also no longer a coexistence of attractors
just to the right of the period-doubling bifurcations.
This suggests that the coexistence
of attractors is less common at lower values of $\epsilon$,
where impact events rob the block of a greater fraction of its energy.

\section{Grazing bifurcations of weakly damped oscillators}
\label{sec:lowDamping}

The previous section described attractors of the
nondimensionalised impact oscillator model \eqref{eq:modelNondim}
with $\zeta = 0$ and $A = A_{\rm graz}$.
These attractors are broadly robust,
so we expect them to persist with $\zeta \approx 0$ and $A \approx A_{\rm graz}$.
Certainly the stable periodic solutions typically
correspond to hyperbolic fixed points of some iterate of the stroboscopic map,
thus vary continuously with respect to $\zeta$ and $A$
by the implicit function theorem \cite{Me17}.
The chaotic attractors described above as being robust to variations in $\omega$,
should also be robust to variations in $\zeta$ and $A$
because the stroboscopic map retains the square-root singularity.
For the fragile chaotic attractors created via period-doubling cascades,
we expect that the entire bifurcation structure shifts with $\zeta$ and $A$,
as long as the overall geometry of the action of the map is unchanged.

\begin{figure}[b!]
\begin{center}
\includegraphics[width=\textwidth]{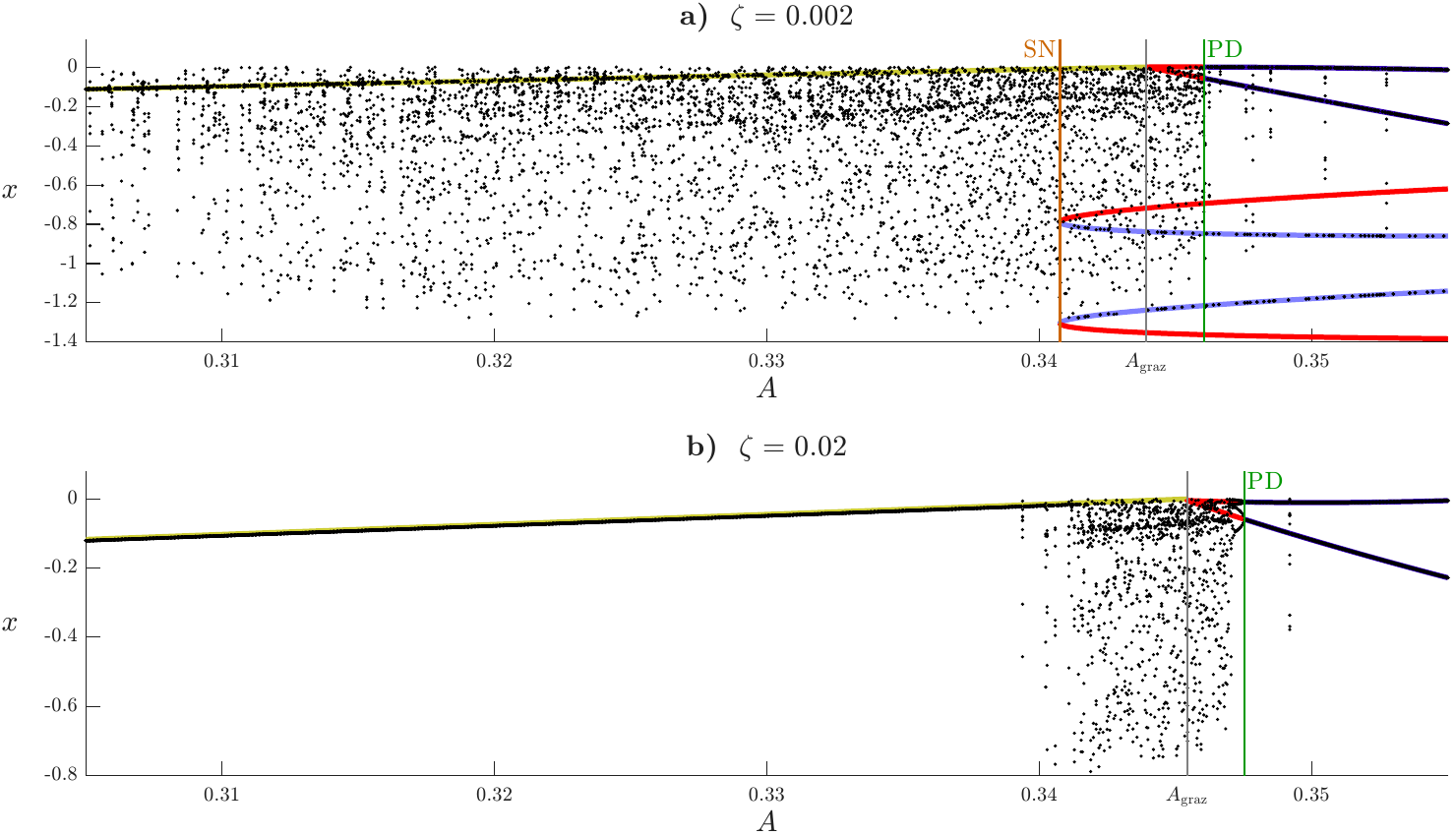}
\caption{
Bifurcation diagrams of the nondimensionalised system \eqref{eq:modelNondim} with
$\omega = 0.81$, $\epsilon = 0.8$, and two different non-zero values of
the damping ratio $\zeta$.
\label{fig:bifDiagThruGrazing}
} 
\end{center}
\end{figure}

As an example, we fix $\omega = 0.81$ and $\epsilon = 0.8$.
These values are similar to those of the physical experiments of
Pavlovskaia {\em et al.}~\cite{PaIn10} that gave a stable period-$\frac{4 \pi}{\omega}$ solution shortly after grazing.
Fig.~\ref{fig:bifDiagThruGrazing}(a) shows how the attractors
of \eqref{eq:modelNondim} vary as we pass through the grazing bifurcation
when the damping ratio is $\zeta = 0.002$.
The yellow line (mostly hidden under the
stroboscopic points of the Monte Carlo simulation)
corresponds the non-impacting solution.

To understand this bifurcation diagram,
recall from Fig.~\ref{fig:bifDiagPPe} that with
$\omega = 0.81$, $\epsilon = 0.8$, $\zeta = 0$, and $A = A_{\rm graz}$,
the system \eqref{eq:modelNondim} has 
a stable two-impact period-$\frac{4 \pi}{\omega}$ solution
coexisting with an apparently chaotic attractor.
And indeed in Fig.~\ref{fig:bifDiagThruGrazing}(a)
at $A = A_{\rm graz}$ the system with $\zeta = 0.002$
has these two attractors, which is consistent with our above comments on robustness.

As the value of $A$ is varied from $A_{\rm graz}$,
these attractors undergo bifurcations and are destroyed.
For instance, the stable period-$\frac{4 \pi}{\omega}$
solution is destroyed in a saddle-node bifurcation,
while the chaotic attractor is destroyed
at $A \approx 0.27$ (outside the plotted range).
Notice also that the grazing bifurcation generates for $A > A_{\rm graz}$
an unstable single-impact period-$\frac{4 \pi}{\omega}$ solution
that quickly gains stability in a period-doubling bifurcation.

With the larger damping ratio $\zeta = 0.02$, see Fig.~\ref{fig:bifDiagThruGrazing}(b),
the two-impact period-$\frac{4 \pi}{\omega}$ solution is no longer present
in the bifurcation diagram.
The chaotic attractor is destroyed at a substantially larger value of $A$ (around $0.34$),
but the behaviour of the single-impact period-$\frac{4 \pi}{\omega}$ solution
is roughly the same.
As in the physical experiments of
Pavlovskaia {\em et al.}~\cite{PaIn10}, the system has a stable period-$\frac{4 \pi}{\omega}$
solution just after the grazing bifurcation.

\section{Discussion}
\label{sec:conc}

We have studied the linear impact oscillator model \eqref{eq:modelOriginal} at grazing in the absence of damping
and discovered through numerical computations a recurring bifurcation structure.
From resonant frequencies there emerge stable periodic solutions
that undergo grazing, saddle-node, and period-doubling bifurcations.
The period-doubling bifurcation initiates a full period-doubling cascade leading a fragile chaotic attractor.
This attractor becomes robust, then contracts to a periodic solution at the next resonant frequency.

We propose that this bifurcation sequence
is the backbone of the overall bifurcation structure of the linear impact oscillator.
Through the examples in \S\ref{sec:lowDamping},
we have shown how the bifurcation structure may differ in the presence of weak damping,
A more thorough investigation for how the arrangement of the bifurcations
changes with damping remains for future work.

Both the stable periodic solution and the robust chaotic attractor that
are created as the value of $\omega$ is varied from a resonant frequency are local objects.
That is, they reside in a neighbourhood of the non-impacting periodic solution.
For this reason, we expect these solutions to arise in the same manner for nonlinear oscillators,
such as those studied in \cite{Fa23,KuNa16}.
Furthermore, we expect that their existence could be proved via a local asymptotic analysis.

Also it remains to understand how the chaos transitions from fragile to robust.
Is there a single bifurcation beyond which periodic windows are no longer dense?
Is the transition driven by global phenomena, akin to the role played by heteroclinic bifurcations
for robust chaos in Lorenz-type systems \cite{DoKr15},
or does the square-root singularity dominate and prevent periodic solutions from being stable?

\section*{Acknowledgements}

This work was supported by Marsden Fund contract MAU2209 managed by Royal Society Te Ap\={a}rangi.

{\footnotesize
\bibliographystyle{unsrt}
\bibliography{NoDamping_BIB}

@book{BlCz99,
	author = {Blazejczyk-Okolewska, B. and Czolczynski, K. and
		Kapitaniak, T. and Wojewoda, J.},
	title = {Chaotic Mechanics in Systems with Impacts and Friction},
	publisher = {World Scientific},
	address = {Singapore},
	year = 1999,
}

@book{Br99,
	author = {Brogliato, B.},
	title = {Nonsmooth Mechanics: Models, Dynamics and Control.},
	publisher = {Springer-Verlag},
	address = {New York},
	year = 1999,
}

@article{ChOt94,
	author = {Chin, W. and Ott, E. and Nusse, H.E. and Grebogi, C.},
	title = {Grazing bifurcations in impact oscillators},
	journal = {Phys. Rev. E},
	volume = 50,
	number = 6,
	pages = {4427-4450},
	year = 1994,
}

@article{DaZh05,
	author = {Dankowicz, H. and Zhao, X.},
	title = {Local analysis of co-dimension-one and co-dimension-two
		grazing bifurcations in impact microactuators},
	journal = {Phys. D},
	volume = 202,
	pages = {238-257},
	year = 2005,
}

@article{DaZh07,
	author = {Dankowicz, H. and Zhao, X. and Misra, S.},
	title = {Near-grazing dynamics in tapping-mode atomic-force microscopy.},
	journal = {Int. J. Non-Linear Mech.},
	volume = 42,
	number = 4,
	pages = {697-709},
	year = 2007,
}

@article{DeVa00,
	author = {de Weger, J. and van de Water, W. and Molenaar, J.},
	title = {Grazing impact oscillations.},
	journal = {Phys. Rev. E},
	volume = 62,
	number = 2,
	pages = {2030-2041},
	year = 2000,
}

@book{DiBu08,
	author = {di Bernardo, M. and Budd, C.J. and Champneys, A.R.
		and Kowalczyk, P.},
	title = {Piecewise-smooth Dynamical Systems. Theory and Applications.},
	publisher = {Springer-Verlag},
	address = {New York},
	year = 2008,
}

@article{DoKr15,
	author = {Doedel, E.J. and Krauskopf, B. and Osinga, H.M.},
	title = {Global organisation of phase space in the transition to
		chaos in the {L}orenz system.},
	journal = {Nonlinearity},
	volume = 28,
	pages = {R113-R139},
	year = 2015,
}

@article{Fa23,
	author = {Farid, M.},
	title = {Dynamics of a hybrid cubic vibro-impact oscillator and nonlinear
		energy sink.},
	journal = {Commun. Nonlinear Sci. Numer. Simul.},
	volume = 117,
	pages = 106978,
	year = 2023,
}

@article{Fo94,
	author = {Foale, S.},
	title = {Analytical determination of bifurcations in an impact oscillator.},
	journal = {Proc. R. Soc. Lond. A},
	volume = 347,
	pages = {353-364},
	year = 1994,
}

@article{HaWi07,
	author = {Halse, C.K. and Wilson, R.E. and di Bernardo, M.
		and Homer, M.E.},
	title = {Coexisting solutions and bifurcations in mechanical
		oscillations with backlash.},
	journal = {J. Sound Vib.},
	volume = 305,
	pages = {854-885},
	year = 2007,
}

@book{Ib09,
	author = {Ibrahim, R.A.},
	title = {Vibro-Impact Dynamics.},
	series = {Lecture Notes in Applied and Computational Mechanics.},
	volume = 43,
	publisher = {Springer},
	address = {New York},
	year = 2009,
}

@article{Iv93,
	author = {Ivanov, A.P.},
	title = {Stabilization Of An Impact Oscillator Near Grazing Incidence
		Owing To Resonance.},
	journal = {J. Sound Vib.},
	volume = 162,
	number = 3,
	pages = {562-565},
	year = 1993,
}

@article{KuNa16,
	author = {Kumar, P. and Narayanan, S. and Gupta, S.},
	title = {Stochastic bifurcations in a vibro-impact
		{D}uffing-{V}an der {P}ol oscillator.},
	journal = {Nonlin. Dyn.},
	volume = 85,
	pages = {439-452},
	year = 2016,
}

@article{LiCh17,
	author = {Li, D. and Chen, H. and Xie, J. and Zhang, J.},
	title = {{Sinai-Ruelle-Bowen} measure for normal form map of grazing
		bifurcations of impact oscillators.},
	journal = {J. Phys. A},
	volume = 50,
	pages = {385103},
	year = 2017,
}

@book{Me17,
	author = {Meiss, J.D.},
	title = {Differential Dynamical Systems.},
	publisher = {SIAM},
	address = {Philadelphia},
	edition = {2nd},
	year = 2017,
}

@article{MiLi19,
	author = {Miao, P. and Li, D. and Yue, Y. and Xie, J. and Grebogi, C.},
	title = {Chaotic attractor of the normal form map for grazing bifurcations
		of impact oscillators.},
	journal = {Phys. D},
	volume = 398,
	pages = {164-170},
	year = 2019,
}

@article{MoCh20,
	author = {Mora, K. and Champneys, A. and Shaw, A. and Friswell, M.I.},
	title = {Explanation of the onset of bouncing cycles in isotropic rotor dynamics;
		a grazing bifurcation analysis.},
	journal = {Proc. R. Soc. A},
	volume = 476,
	pages = {20190549},
	year = 2020,
}

@article{No91,
	author = {Nordmark, A.B.},	
	title = {Non-periodic motion caused by grazing incidence in impact oscillators.},
	journal = {J. Sound Vib.},
	volume = 145,
	number = 2,
	pages = {279-297},
	year = 1991,
}

@article{No97,
	author = {Nordmark, A.B.},
	title = {Universal limit mapping in grazing bifurcations.},
	journal = {Phys. Rev. E},
	volume = 55,
	number = 1,
	pages = {266-270},
	year = 1997,
}

@article{No01,
	author = {Nordmark, A.B.},
	title = {Existence of periodic orbits in grazing bifurcations
		of impacting mechanical oscillators.},
	journal = {Nonlinearity},
	volume = 14,
	pages = {1517-1542},
	year = 2001,
}

@article{PaIn10,
	author = {Pavlovskaia, E. and Ing, J. and Wiercigroch, M. and
		Banerjee, S.},
	title = {Complex Dynamics of Bilinear Oscillator Close to Grazing.},
	journal = {Int. J. Bifurcation Chaos},
	volume = 20,
	number = 11,
	pages = {3801-3817},
	year = 2010,
}

@article{Pe96,
	author = {Peterka, F.},
	title = {Bifurcations and Transition Phenomena in an Impact Oscillator.},
	journal = {Chaos Solitons Fractals},
	volume = 7,
	number = 10,
	pages = {1635-1647},
	year = 1996,
}

@article{PeVa92,
	author = {Peterka, F. and Vacik, J.},
	title = {Transition to Chaotic Motion in Mechanical
		Systems with Impacts.},
 	journal = {J. Sound Vib.},
	volume = 154,
	number = 1,
	pages = {95-115},
	year = 1992,
}

@unpublished{SiGh26,
	author = {Simpson, D.J.W. and Ghosh, I.},
	title = {Resonant grazing bifurcations revisited.},
	year = 2026,
	note = {{\em https://arxiv.org/abs/2602.22797}},
}

@article{ThZh06,
	author = {Thota, P. and Zhao, X. and Dankowicz, H.},
	title = {Co-dimension-two Grazing Bifurcations in Single-degree-of-freedom
		Impact Oscillators.},
	journal = {J. Comput. Nonlinear Dynam.},
	volume = 1,
	number = 4,
	pages = {328-335},
	year = 2006,
}

@article{WiVi14,
	author = {Witelski, T. and Virgin, L.N. and George, C.},
	title = {A driven system of impacting pendulums:
		{E}xperiments and simulations.},
	journal = {J. Sound Vib.},
	volume = 333,
	pages = {1734-1753},
	year = 2014,
}
}

\end{document}